\documentclass[10pt]{article}
\usepackage{amsmath}
\usepackage{amsthm}
\usepackage{amssymb}
\usepackage{amscd}
\usepackage{diagrams}
\theoremstyle{plain}

\newtheorem{thm}[subsection]{Theorem}
\newtheorem{prp}[subsection]{Proposition}
\newtheorem{cor}[subsection]{Corollary}
\newtheorem{lma}[subsection]{Lemma}
\theoremstyle{remark}
\newtheorem{rmk}[subsection]{Remark}
\newtheorem{exms}[subsection]{Examples}
\newtheorem{tit}[subsection]{}
\newtheorem{stit}[subsubsection]{}
\def\EE{\mathcal{E}}

\def\TT{\mathbb{T}}
\def\NN{\mathrm{N}}
\def\CC{\mathcal{C}}
\def\DD{\mathcal{D}}
\def\EE{\mathcal{E}}
\def\KK{\mathcal{K}}
\def\LL{\mathcal{L}}

\def\UU{\mathcal{U}}

\def\FF{\mathcal{F}}

\def\AA{\mathcal{A}}

\def\ZZ{\mathbb{Z}}
\def\WW{\mathrm{W}}

\def\Coll{\mathrm{Coll}}
\def\Oper{\mathrm{Oper}}
\def\Hopf{\mathrm{Hopf}}
\def\Path{\mathrm{Path}}
\def\Sh{\mathrm{Sh}}
\def\Ch{\mathrm{Ch}}
\def\Sets{\mathrm{Sets}}
\def\Com{\mathcal{C}om}
\def\Ass{\mathcal{A}ss}
\def\op{\mathrm{op}}

\def\Hom{\mathrm{Hom}}

\def\Alg{\mathrm{Alg}}
\def\Coalg{\mathrm{Coalg}}
\def\Ho{\mathbf{Ho}}
\def\Sg{\Sigma}
\def\sg{\sigma}
\def\ito{\rightarrowtail}
\def\cto{\twoheadrightarrow}

\def\lrto{\leftrightarrows}

\def\eqv{\overset{\sim}{\longrightarrow}}

\def\inp{\mathit{in}}
\def\Aut{\mathrm{Aut}}

\begin{document}
\title{Axiomatic homotopy theory for operads}
\author{Clemens Berger and Ieke Moerdijk}
\date{}
\maketitle

\begin{abstract}We give sufficient conditions for the existence of a model structure on operads in an arbitrary symmetric monoidal model category. General invariance properties for homotopy algebras over operads are deduced.\end{abstract}

\section*{Introduction}

Operads are a device for describing algebraic structures. Initially, they served to define algebraic structures on topological spaces with constraints holding only up to a coherent system of homotopies. Stasheff's $A_\infty$-spaces resp. Boardman, Vogt and May's $E_\infty$-spaces are spaces with a ``homotopy associative'' resp. ``homotopy commutative'' multiplication, cf. \cite{St}, \cite{BV}, \cite{M}. The corresponding $A_\infty$- resp. $E_\infty$-operad is a deformation of the operad acting on strictly associative resp. strictly commutative monoids. This method of deforming algebraic structures via operads has recently received new interest, cf. Mandell \cite{Man}, Kontsevich-Soibelman \cite{KS}, Hinich \cite{H}, and others.

Symmetric operads may be defined in any symmetric monoidal category. We propose here a general homotopy theory for such operads by means of Quillen closed model categories \cite{Q}. We show that under certain conditions the (reduced) operads of a symmetric monoidal model category carry a model structure, with weak equivalences and fibrations defined on the level of the underlying collections. These conditions concern the existence of a suitable ``interval'' with comultiplication; they are easy to verify, and are well known to hold in many standard situations. In particular, they hold for topological, simplicial and chain operads. Our approach may be compared with Hinich \cite{H}, Spitzweck \cite{Sp} and Vogt \cite{V}, but is much more elementary, since it relies on the existence of path-objects rather than on an intricate analysis of pushouts. Our method also immediately extends to coloured operads.

The principal benefit of a model structure on operads is an intrinsic definition of homotopy algebras over an operad, namely as the algebras over a cofibrant replacement of the given operad, cf. Markl \cite{Mkl}. The algebras over cofibrant operads carry a model structure for which a variant of the Boardman-Vogt homotopy invariance property holds. A larger class of operads, here called admissible $\Sigma$-cofibrant, will be shown to induce the same homotopy theory for their algebras as any of their cofibrant replacements. This is important, since most of the commonly used $E_\infty$-operads are actually $\Sigma$-cofibrant, but not cofibrant. As main comparison theorem, we show that the base-change adjunction with respect to a weak equivalence of admissible $\Sigma$-cofibrant operads induces an equivalence of the corresponding homotopy categories of algebras. We also compare homotopy algebras in different symmetric monoidal model categories: for instance, the homotopy category of simplicial $E_\infty$-algebras is equivalent to the homotopy category of topological $E_\infty$-algebras.\vspace{3ex}

The plan of this article is as follows:\vspace{1ex}

Section $1$ first reviews the basic concepts involving operads and algebras over an operad in an arbitrary closed symmetric monoidal category. We then discuss the dual notions of a cooperad and a coalgebra and define two convolution pairings subsequently used for the construction of path-objects.

Section $2$ recalls the basics of (monoidal) model categories with special emphasis on the transfer of model structures.

Section $3$ establishes the two main theorems giving sufficient conditions for the existence of a model structure on operads with weak equivalences and fibrations defined at the level of the underlying collections. We also discuss Boardman and Vogt's $\WW$-construction as well as a model-theoretic formulation of their homotopy invariance property. We finally discuss the standard examples (simplicial, topological, chain and sheaf operads) where our method yields model structures.

Section $4$ contains two comparison theorems: the first shows that the base-change adjunction with respect to a weak equivalence of admissible $\Sg$-cofibrant operads is a Quillen equivalence. The second shows that under mild assumptions, a monoidal Quillen equivalence between monoidal model categories induces a Quillen equivalence between mutually corresponding categories of ``homotopy algebras''.

The Appendix contains complete proofs for some key properties of $\Sg$-cofibrant operads used in Section $4$.\vspace{3ex}

{\sc Acknowledgements:} The research for this paper started during the Workshop on Operads in Osnabr\" uck (December 2000), and we thank R. Vogt for inviting us both. In addition, we would like to thank B. Fresse, P. van der Laan, M. Markl, P. May, B. Shipley, J. Smith and B. Toen for helpful discussions. The detailed comments of the referee have been much appreciated. This work has been supported by the SFB at Bielefeld/Osnabr\"uck, the EU-project Modern Homotopy Theory and the Netherlands Science Organisation (NWO).\newpage

\section{Operads, algebras and convolution products}

The main purpose of this section is to discuss the ``convolution operad'' associated to a cooperad and an operad, as well as the ``convolution algebra'' associated to a coalgebra and an algebra. We begin by recalling some standard notation and terminology concerning operads and algebras.

Throughout this paper, $\EE=(\EE,\otimes,I,\tau)$ is a fixed closed symmetric monoidal category. We assume that $\EE$ has small colimits and finite limits. The closedness of $\EE$ means that the functor $-\otimes X$ has a right adjoint, denoted $(-)^X$. Every symmetric monoidal category is equivalent to one with a strictly associative monoidal structure. Therefore, the bracketing convention for multiple tensor products is not really important; we adopt a bracketing from left to right, which seems best suited with respect to the adjunction with the internal hom.

For a discrete group $G$, we write $\EE^G$ for the category of objects in $\EE$ with a right $G$-action. It is again a closed symmetric monoidal category, and the forgetful functor $\EE^G\to\EE$ preserves this structure and has a left adjoint, denoted $(-)[G]$. This applies in particular to each of the symmetric groups $\Sg_n$, where for consistency $\Sg_0$ and $\Sg_1$ both denote the trivial group. The product of the categories $\EE^{\Sg_n}$ is called the category of \emph{collections}, and denoted$$\Coll(\EE)=\prod_{n\geq 0}\EE^{\Sg_n}.$$Its objects are written $P=(P(n))_{n\geq 0}$. Each collection $P$ induces an endofunctor (again denoted) $P:\EE\to\EE$, by$$P(X)=\coprod_{n\geq 0}P(n)\otimes_{\Sg_n}X^{\otimes n}.$$This endofunctor has the structure of a \emph{monad} if the defining collection is an \emph{operad}, which means that $P$ comes equipped with a unit $I\to P(1)$ and with a family of structure maps$$P(k)\otimes P(n_1)\otimes\cdots\otimes P(n_k)\to P(n_1+\cdots+n_k)$$satisfying well known equivariance, associativity and unit conditions. For more details, see e.g. Boardman-Vogt \cite[lemma 2.43]{BV}, May \cite[def. 1.1]{M} and Getzler-Jones \cite[def. 1.4]{GJ}. The category of operads in $\EE$ is denoted by $\Oper(\EE)$. 

A \emph{cooperad} is a collection $C$ equipped with a counit $C(1)\to I$ and structure maps $C(n_1+\cdots+ n_k)\to C(k)\otimes C(n_1)\otimes\cdots\otimes C(n_k)$ satisfying the dual conditions. 

If $C$ is a cooperad and $P$ is an operad, then the collection $P^C$ defined by \begin{gather*}\label{convolution}P^C(n)=P(n)^{C(n)}\end{gather*}(with the usual $\Sg_n$-actions by conjugation on the exponent) has a natural convolution operad structure with structure maps given by\begin{align*}P^C(k)\otimes P^C(n_1)\otimes\cdots\otimes P^C(n_k)&\to (P(k)\otimes\cdots\otimes P(n_k))^{(C(k)\otimes \cdots\otimes C(n_k))}\\&\to P(n_1+\cdots+n_k)^{C(n_1+\cdots+n_k)}\\&=P^C(n_1+\cdots+n_k)\end{align*}

A (commutative) \emph{Hopf object} is an object $H=(H,m,\eta,\Delta,\epsilon)$  such that $(H, m, \eta)$ is a (commutative) monoid, $(H,\Delta, \epsilon)$ a comonoid and $m, \eta$  are maps of comonoids (resp. $\Delta, \epsilon$ maps of monoids). Here, the symmetry $\tau$ of $\EE$ enters in an essential manner. The category of commutative Hopf objects in $\EE$ is denoted by $\Hopf(\EE)$. If the tensor of $\EE$ is the cartesian product, any monoid $H$ has a canonical Hopf structure, given by the diagonal $\Delta:H\to H\times H$ and the unique map $\epsilon:H\to I$ to the terminal object. 

Each commutative Hopf object $H$ defines a \emph{cooperad} $TH$ with underlying collection given by $(TH)(n)=H^{\otimes n}$. For $n=n_1+\cdots+n_k$, the structure map\begin{gather*}\label{cooperad}H^{\otimes n }\to H^{\otimes k}\otimes H^{\otimes n_1}\otimes\cdots\otimes H^{\otimes n_k}\end{gather*} is the composite of the comonoid structure $H^{\otimes n}\to H^{\otimes n}\otimes H^{\otimes n}$ on $H^{\otimes n}$ with the map $p\otimes i$, where $i:H^{\otimes n}\cong H^{\otimes n_1}\otimes\cdots\otimes H^{\otimes n_k}$ is the canonical isomorphism, and $p: H^{\otimes n}\to H^{\otimes k}$ is the composite of $i$ with the product on each of the $k$ tensor factors. The commutativity of $H$ guarantees that the equivariance conditions for the cooperad $TH$ are satisfied.

For any operad $P$, a \emph{$P$-algebra} $A$ is an object of $\EE$ equipped with structure maps $P(n)\otimes A^{\otimes n}\to A,\, n\geq 0,$ satisfying well known equivariance, associativity and unit conditions, cf. \cite{BV}, \cite{M}, \cite{GJ}. Dually, a \emph{$P$-coalgebra} $B$ is an object of $\EE$ equipped with structure maps $P(n)\otimes B\to  B^{\otimes n}$ satisfying the dual conditions.

We denote the category of $P$-algebras by $\Alg_P$ and the category of $P$-coalgebras by $\Coalg_P$. A $P$-algebra structure on $A$ corresponds also to an operad map $P\to\EE_A$ with values in the \emph{endomorphism-operad}, defined by $\EE_A(n)=A^{(A^{\otimes n})}$ with the natural compositional operad structure, cf. Smirnov \cite{Sm}. Dually, a $P$-coalgebra structure on $B$ corresponds to an operad map $P\to\EE_B^\op$ with values in the \emph{coendomorphism-operad}, defined by $\EE^\op_B(n)=(B^{\otimes n})^B$. The product $P\otimes Q$ of two operads is defined by $(P\otimes Q)(n)=P(n)\otimes Q(n)$ with the obvious structure maps.

\begin{prp}\label{Hopf}There are natural convolution pairings\begin{gather*}\Hopf(\EE)^\op\times\Oper(\EE)\to\Oper(\EE)\\\Coalg_Q^\op\times\Alg_P\to\Alg_{P\otimes Q}\end{gather*}\end{prp}
\begin{proof} The first pairing maps a commutative Hopf object $H$ and an operad $P$ to the convolution operad $P^{TH}$. The second pairing maps a $Q$-coalgebra $B$ and a $P$-algebra $A$ to the object $A^B$ equipped with the following $P\otimes Q$-algebra structure:\begin{align*}(P\otimes Q)(n)\otimes (A^B)^{\otimes n}&\to P(n)\otimes (A^B)^{\otimes n}\otimes Q(n)\\&\to A^{(A^{\otimes n})}\otimes(A^{\otimes n})^{(B^{\otimes n})}\otimes(B^{\otimes n})^B\\&\to A^B\end{align*}\end{proof}

Note that, in particular, if $P$ is an operad with diagonal $\Delta:P\to P\otimes P$ (e.g. $P$ is a Hopf operad), the second pairing for $P=Q$ together with the pullback along $\Delta$ defines a convolution product $\Coalg_P^\op\times\Alg_P\to\Alg_P$.

There is also a convolution product of a $P$-algebra $A$ and a $C$-coalgebra $B$ (for a cooperad $C$) yielding a $P^C$-algebra $A^B$, but we will not use this construction in this paper.

\section{Background on model categories}

\begin{tit}\emph{Model categories.}\label{model}\end{tit}

In this paper, a \emph{model category} always means what Quillen calls a closed model category \cite[I.5]{Q}. An adjoint pair between model categories is a \emph{Quillen pair} if the left adjoint preserves cofibrations and the right adjoint preserves fibrations. This implies that the adjoint pair passes to the homotopy categories. If the derived adjoint pair between the homotopy categories is an equivalence, the original Quillen pair is called a \emph{Quillen equivalence}.

A model category is \emph{left proper}, if the class of weak equivalences is closed under pushouts along cofibrations. A sufficient condition for left properness is that all objects of the model category be cofibrant.

In any model category, the colimit functor sends weak equivalences of directed Reedy-cofibrant diagrams to weak equivalences, cf. \cite[prp. 11.5]{DHK}. This implies in particular that a horizontal ladder of cofibrations between cofibrant objects in which all vertical maps are weak equivalences yields a weak equivalence in the colimit, a fact needed in the Appendix.

\begin{tit}\emph{Monoidal model categories.}\label{monoid}\end{tit}

A \emph{monoidal model category} is a closed symmetric monoidal category endowed with a model structure subject to the following \emph{pushout-product axiom} (cf. \cite{SS}):  

For any pair of cofibrations $f:X\ito Y$ and $f':X'\ito Y'$, the induced map $(X\otimes Y')\cup_{X\otimes X'} (Y\otimes X')\to Y\otimes Y'$ is a cofibration, which is trivial if $f$ or $f'$ is trivial. In particular, tensoring with cofibrant objects preserves cofibrations and trivial cofibrations. However, the tensor product of two (trivial) cofibrations is in general not a (trivial) cofibration, cf. (\ref{stable}).\vspace{1ex} 


We shall repeatedly use the following basic lemma:

\begin{lma}\label{weak}Let $f:X\to Y$ be a map between cofibrant objects of a monoidal model category. If $f$ is a weak equivalence, then for every fibrant object $Z$, the induced map $f^*:Z^Y\to Z^X$ is a weak equivalence. The converse holds as soon as the unit of the monoidal model category is cofibrant.\end{lma}

\begin{proof} For a trivial cofibration $f$, the first assertion is true by exponential transposition and the dual of pushout-product axiom. K. Brown's lemma allows us to conclude the proof, cf. \cite[lemma 1.1.12]{Hov}. Conversely, assume that $f^*$ is a weak equivalence for fibrant objects $Z$. Then, since the unit $I$ of the monoidal model category $\EE$ is cofibrant, $f^*$ induces a bijection $\Ho(\EE)(I,Z^Y)\to\Ho(\EE)(I,Z^X)$ in the homotopy category. The tensor-hom adjunction is compatible with the homotopy relation, so that we obtain a bijection $\Ho(\EE)(Y,Z)\to\Ho(\EE)(X,Z)$ for every fibrant object $Z$. Therefore, $\Ho(\EE)(f)$ is an isomorphism, which shows that $f$ is a weak equivalence.\end{proof}

\begin{tit}\emph{Symmetric monoidal functors and fibrant replacement functors.}\label{monoidal}\end{tit}

A functor $F:(\DD,\otimes_\DD,I_\DD,\tau_\DD)\to(\EE,\otimes_\EE,I_\EE,\tau_\EE)$ between symmetric monoidal categories is \emph{symmetric monoidal} if $F$ comes equipped with a unit $I_\EE\to F(I_\DD)$ and a binatural transformation $F(X)\otimes_\EE F(Y)\to F(X\otimes_\DD Y)$ satisfying familiar unit, associativity and symmetry conditions, cf. \cite[III.20]{MMSS}. A symmetric monoidal functor maps commutative monoids to commutative monoids, and operads to operads. A symmetric monoidal functor is said to be \emph{strong} if the structure maps are isomorphisms. For instance, any product-preserving functor between cartesian closed categories is strong symmetric monoidal.

A \emph{fibrant replacement} for an object $X$ is a weak equivalence $X\eqv \tilde{X}$ with fibrant codomain. If this weak equivalence is part of a natural transformation $id_\EE\to(-)^{\tilde{}}$ we say that the model category admits a \emph{fibrant replacement functor}. This fibrant replacement functor is called symmetric monoidal if the functor $(-)^{\tilde{}}$ is symmetric monoidal and the structure maps $\tilde{X}\otimes\tilde{Y}\to(X\otimes Y)^{\tilde{}}$ are defined under $X\otimes Y$.

\begin{tit}\emph{Cofibrant generation and transfer of model structures.}\label{transfer}\end{tit}
A model category is \emph{cofibrantly generated} if the category is cocomplete and admits generating \emph{sets} of (trivial) cofibrations with small domains, cf. Dwyer-Hirschhorn-Kan \cite[II.7.4]{DHK} and Hovey \cite[II.1]{Hov}, who use a slightly more general concept. ``Small'' means $\lambda$-small for some regular cardinal $\lambda$ and ``generating'' means that the fibrations (resp. trivial fibrations) are characterized by their right lifting property with respect to the generating trivial cofibrations (resp. cofibrations).

With respect to a given set of generating (trivial) cofibrations, a \emph{(trivial) cellular extension}  is a sequential colimit of pushouts of generating (trivial) cofibrations. A \emph{cellular object} is a cellular extension of the initial object. In a cofibrantly generated model category, each (trivial) cofibration is a codomain-retract of a (trivial) cellular extension. In particular, each cofibrant object is a retract of a cellular object.

Cofibrantly generated model structures may be transferred along the left adjoint functor of an adjunction. The first general statement of such a transfer in the literature is due to Crans \cite{C}.\vspace{1ex}  

\emph{Transfer principle}: Let $\DD$ be a cofibrantly generated model category and let $F:\DD\lrto\EE:G$ be an adjunction with left adjoint $F$ and right adjoint $G$. Assume that $\EE$ has small colimits and finite limits. Define a map $f$ in $\EE$ to be a weak equivalence (resp. fibration)  iff $G(f)$ is a weak equivalence (resp. fibration). Then this defines a cofibrantly generated model structure on $\EE$ provided\vspace{1ex}

(i) the functor $F$ preserves small objects;

(ii) any sequential colimit of pushouts of images under $F$ of the generating trivial cofibrations of $\DD$ yields a weak equivalence in $\EE$.\vspace{1ex}

 Condition (i) holds for instance if $G$ preserves filtered colimits.

\begin{tit}\emph{Quillen's path-object argument.}\label{path}\end{tit}

In practice, the condition (ii) above is the crucial property to be verified. This is often hard, but there is one special case in which an argument of Quillen's yields (ii). Recall that a \emph{path-object} for $X$ is a factorisation of its diagonal into a weak equivalence followed by a fibration $X\eqv\Path(X)\cto X\times X$. Assume, the following two conditions hold:\vspace{1ex}

(a) $\EE$ has a fibrant replacement functor;

(b) $\EE$ has functorial path-objects for fibrant objects.\vspace{1ex}

Then condition (ii) for transfer is satisfied, cf. \cite[II.4]{Q}, \cite[7.6]{R}, \cite[A.3]{SS}. If all objects in $\EE$ are fibrant, (a) of course becomes redundant, since the identity serves as fibrant replacement functor.\vspace{1ex}

\section{Model structure on operads}

This section gives sufficient conditions for the category of operads of a monoidal model category to have a model structure. We shall see below that these conditions are easy to verify in many standard examples. If the model structure on $\EE$ is cofibrantly generated, the transfer principle (\ref{transfer}) implies that for any discrete group $G$, the model structure on $\EE$ lifts to a model structure on $\EE^G$, in which a map is a weak equivalence (resp. fibration) iff it is so once we forget the $G$-action. It follows that $\Coll(\EE)$ is a cofibrantly generated model category, in which a map $P\to Q$ is a weak equivalence (resp. fibration) iff for each $n$, the map $P(n)\to Q(n)$ is a weak equivalence (resp. fibration) in $\EE$. Using the \emph{path-object argument} (\ref{path}), we shall transfer this model structure along the free-forgetful adjunction\begin{gather*}\label{adjunction}\FF:\Coll(\EE)\lrto\Oper(\EE):\UU\end{gather*}(or a similar adjunction between reduced collections and reduced operads as defined below).

In the underlying category $\EE$, the unit $I$ is a commutative Hopf object, by the canonical isomorphism $I\otimes I\to I$ and its inverse. The coproduct $I\sqcup I$ is also a Hopf object in a natural way. Indeed, if we label the two summands with the elements of $\ZZ/2\ZZ$, then the multiplication on $I\sqcup I$ is induced by that of $I$ and the one on $\ZZ/2\ZZ$, while the comultiplication is induced by that of $I$ and the diagonal of $\ZZ/2\ZZ$. Furthermore, the folding map $I\sqcup I\to I$ is a map of Hopf objects. We say that $\EE$ admits a (commutative) \emph{Hopf interval} if this folding map can be factored into a cofibration followed by a weak equivalence\begin{gather*}\label{interval}I\sqcup I\ito H\eqv I\end{gather*}where $H$ is a (commutative) Hopf object and both maps are maps of Hopf objects.\vspace{1ex}

 The homotopy theory of operads suffers from the fact that by the very definition of an operad, the $0$-th term $P(0)$ is the initial object of the category of $P$-algebras, and that moreover the $P$-algebras under a fixed $P$-algebra $A$ are the algebras for another operad, whose $0$-th term is $A$. Therefore, the homotopy theory of operads \emph{subsumes} the homotopy theory of algebras over a given operad, and inherits the difficulties of the latter. For instance, for monoidal model categories which are not cartesian closed, commutative monoids in general do not admit a well behaved homotopy theory, so that general operads do not have a well behaved homotopy theory either. In order to avoid this mixture of the operad and algebra levels, we introduce the category of \emph{reduced} operads. An operad $P$ is reduced if $P(0)$ is the unit of $\EE$. A map of reduced operads is a map of operads $\phi:P\to Q$ such that $\phi(0):P(0)\to Q(0)$ is the identity. Observe that the action (\ref{Hopf}) of commutative Hopf objects on operads restricts to reduced operads, and that the collection underlying a reduced operad is in fact a collection in $\EE/I$.

\begin{thm}\label{operad}Let $\EE$ be a monoidal model category with unit $I$ such that\begin{itemize}\item $\EE$ is cofibrantly generated and its unit is cofibrant;\item $\EE/I$ has a symmetric monoidal fibrant replacement functor;\item $\EE$ admits a commutative Hopf interval.\end{itemize}
Then, there is a cofibrantly generated model structure on the category of reduced operads, in which a map $P\to Q$ is a weak equivalence (resp. fibration) iff for each $n>0$, the map $P(n)\to Q(n)$ is a weak equivalence (resp. fibration) in $\EE$.\end{thm}

\begin{proof}We shall construct the model structure on operads by transfer (\ref{transfer}) using the path-object argument (\ref{path}). The category of reduced operads is cocomplete and finitely complete, since the same is true for the category of reduced collections in $\EE/I$, and since the forgetful functor from reduced operads to reduced collections in $\EE/I$ is monadic and preserves filtered colimits.

Let $P$ be a reduced operad and $\tilde{P}$ be the collection defined by $\tilde{P}(0)=P(0)$ and $\tilde{P}(n)=P(n)\tilde{}$ if $n>0$, where $X\mapsto X\tilde{}$ is the symmetric monoidal fibrant replacement functor in $\EE/I$. Then the operad structure on $P$ induces an operad structure on $\tilde{P}$, so that $\tilde{P}$ is a fibrant replacement for $P$ in the category of reduced operads. Thus, (\ref{path}a) holds.

Assume now that $P$ is a fibrant reduced operad. The construction (\ref{Hopf}) applied to the Hopf interval $H$ and to $P$ yields a functorial path-object: $$P=P^{TI}\eqv P^{TH}\cto P^{T(I\sqcup I)}\cto P^I\times P^I=P\times P.$$ Indeed, the $n$-fold tensor product $(I\sqcup I)^{\otimes n}\to H^{\otimes n}$ is a cofibration by the pushout-product axiom and the assumption that $I$ is cofibrant. Therefore, $P^{TH}\to P^{T(I\sqcup I)}$ is a fibration. The canonical map $P^{T(I\sqcup I)}\to P^I\times P^I$ induces for each $n$ the projection $P(n)^{(2^n)}\to P(n)^2$ onto the first and last factor. Since $P$ is fibrant, this is a fibration for $n\geq 1$. This shows that $P^{TH}\to P\times P$ is a fibration. Moreover, since $H\eqv I$ admits a trivial cofibration as section, the $n$-fold tensor product $H^{\otimes n}\to I^{\otimes n}$ is a weak equivalence between cofibrant objects, whence $P\to P^{TH}$ is a weak equivalence by (\ref{weak}). Thus, (\ref{path}b) holds and the transfer applies.\end{proof}

If the monoidal structure is cartesian closed, we can weaken the hypotheses. Furthermore, the restriction to \emph{reduced} operads isn't necessary in this case:

\begin{thm}\label{operad2}Let $\EE$ be a cartesian closed model category such that\begin{itemize}\item $\EE$ is cofibrantly generated and the terminal object of $\EE$ is cofibrant;\item $\EE$ has a symmetric monoidal fibrant replacement functor.\end{itemize}
Then, there is a cofibrantly generated model structure on the category of operads, in which a map $P\to Q$ is a weak equivalence (resp. fibration) iff for each $n$, the map $P(n)\to Q(n)$ is a weak equivalence (resp. fibration) in $\EE$.\end{thm}

\begin{proof}The first part is identical to the preceding proof, except that we put $\tilde{P}(0)=P(0)\tilde{}$ this time, and use unreduced operads and collections. For the construction of a functorial path-object for fibrant operads, we use the fact that in a cartesian closed category exponentiation is product-preserving and hence strong symmetric monoidal. This implies that for any interval $I\sqcup I\ito J\eqv I$, mapping into a fibrant operad $P$ yields a path-object in the category of operads: $P=P^I\eqv P^J\cto P^{I\sqcup I}\cong P\times P$.\end{proof}

The above proof can be adapted to a slightly more general situation: Indeed, Theorem (\ref{operad2}) remains true for a cofibrantly generated monoidal model category $\EE$ having a cofibrant unit, a symmetric monoidal fibrant replacement functor and an interval with a coassociative and cocommutative comultiplication.

\begin{exms}\label{examples}--

\begin{stit}\emph{Simplicial operads.}

The category of simplicial sets is a cartesian closed model category, in which the cofibrations are the monomorphisms and the weak equivalences are the realization weak equivalences. The pushout-product axiom is well known in this case.  The model structure is cofibrantly generated by the boundary-inclusions (resp. horn-inclusions) of the standard $n$-simplices. As symmetric monoidal fibrant replacement functor, we can choose either Kan's $Ex^\infty$ functor or the singular complex of the geometric realization functor, since both are product-preserving. Therefore, simplicial operads form a model category by (\ref{operad2}).

This fact could also have been deduced from Quillen's theorem \cite[II.4]{Q} that the simplicial objects of a (cocomplete, finitely complete) category with a set of small projective generators admit a canonical model structure. Indeed, operads in sets form such a category. Since the projective generators represent evaluation at $n$, Quillen's structure coincides with ours. There is \emph{another} model structure on simplicial operads, obtained by restricting Rezk's model structure \cite[7.5]{R} on simplicial theories to simplicial operads. The class of weak equivalences for this model structure is smaller: in particular, two simplicial operads are weakly equivalent for Rezk's model structure if and \emph{only if} they define equivalent homotopy categories of algebras \cite[8.6]{R}, which is in general not true for our model structure, cf. however (\ref{invariant}) below.\end{stit}

\begin{stit}\emph{Topological operads.}

The category of compactly generated spaces is a cartesian closed model category, in which the weak equivalences are the weak homotopy equivalences and the fibrations are the Serre fibrations, cf. \cite[II.3]{Q}. The pushout-product axiom follows from the fact that this model structure on topological spaces is obtained by transfer from the model structure on simplicial sets along a product-preserving realization functor. The model structure is cofibrantly generated by the sphere (resp. ball) inclusions and all objects are fibrant. Therefore, topological operads form a model category by (\ref{operad2}).\end{stit}

\begin{stit}\emph{Chain operads.}

For any commutative ring $R$ with unit, the category $\Ch(R)$ of $\ZZ$-graded chain complexes of $R$-modules is a cofibrantly generated monoidal model category with quasi-isomorphisms as weak equivalences and epimorphisms as fibrations, cf. \cite[2.3]{Hov}. The normalized $R$-chain functor $\NN^R_*:\Sets^{\Delta^{op}}\to\Ch(R)$ is symmetric monoidal. The structure maps are given by the Eilenberg-Zilber shuffle map. The normalized $R$-chain functor also carries a comonoidal structure given by the Alexander-Whitney diagonal; the latter is however not symmetric. Moreover, there is a mixed associativity relation relating the monoidal and the comonoidal structure, which implies that $\NN^R_*$ sends Hopf objects to Hopf objects; in particular, the image $\NN^R_*(\Delta[1])$ of the standard simplicial interval $\Delta[1]$ is a commutative (but not cocommutative) Hopf interval in the category of $R$-chain complexes. The unit of $\Ch(R)$ is $R$ concentrated in degree $0$, which is clearly cofibrant. All objects have a fibrant replacement over the unit. Therefore, reduced $R$-chain operads carry a model structure by (\ref{operad}).

 This has been proved by Hinich using a different method \cite{H}. Hinich is not explicit about the $0$-th term of his operads. There is however no transferred model structure for unreduced $R$-chain operads, because the coproduct with the operad $\Com$ for commutative $R$-chain algebras does not send trivial cofibrations of \emph{unreduced} operads to weak equivalences.\end{stit}

\begin{stit}\emph{Sheaf operads.}

Generalizing the examples above, we can consider simplicial objects (resp. $R$-chain complexes) in the topos $\Sh(\TT)$ of set-valued sheaves on a small site $\TT$ of finite type. According to a widely circulated letter from A. Joyal to A. Grothendieck, there is a cofibrantly generated monoidal model structure on both categories with monomorphisms as cofibrations and ``stalkwise weak equivalences'' as weak equivalences. In the simplicial case, the pushout-product axiom is easy to verify; in the $R$-chain-complex case, the pushout-product axiom only holds if $R$ is a field. Morel-Voevodsky \cite[2.1.66]{MV} show that the category of simplicial sheaves on a site $\TT$ of finite type admits a symmetric monoidal fibrant replacement functor. If $R$ is a field, such a replacement functor can also be constructed for $R$-chain complexes in $\Sh(\TT)$. Moreover, in the latter case, the constant sheaf at $\NN_*^R(\Delta[1])$ is a Hopf interval. Therefore, simplicial operads (resp. reduced $R$-chain operads) in $\Sh(\TT)$ carry a model structure by (\ref{operad2}) (resp. \ref{operad}). In particular, there exists a ``continously varying'' simplicial (resp. $R$-chain) $E_\infty$-operad on $\TT$.\end{stit}\end{exms}\vspace{1ex}

\begin{rmk}Boardman and Vogt \cite[III.1]{BV} define for each topological operad $P$ an operad $\WW P$ with the property that $\WW P$-algebras may be considered as homotopy $P$-algebras. In Section $4$, we define homotopy $P$-algebras as the algebras over a cofibrant replacement for $P$. The relationship between these two notions of homotopy algebras is established by the following property of the $\WW$-construction, where a topological operad $P$ is called \emph{well-pointed} (resp. \emph{$\Sg$-cofibrant}) if the unit $I\to P(1)$ is a cofibration (resp. the underlying collection is cofibrant):

\hspace{3em}\emph{for any well-pointed $\Sg$-cofibrant operad $P$, the augment-} 

\hspace{3em}\emph{ation $\epsilon_P:\WW P\eqv P$ is a cofibrant replacement for $P$.}

\noindent This statement is essentially proved by Vogt \cite[thm. 4.1]{V}. We have shown that the $\WW$-construction for operads can be defined in any monoidal model category $\EE$ equipped with a suitable interval $H$, and that it defines a functorial cofibrant replacement in this general context. The proof is rather technical, and will be presented in \cite{BM}.\end{rmk}
 
Endomorphism-operads are not reduced, since $\EE_X(0)=X$. However, any object $X$ \emph{under} $I$ defines a \emph{reduced} endomorphism-operad $\bar{\EE}_X$. If $P$ is reduced, a $P$-algebra structure on $X$ is also equivalent to a base point $I\to X$ together with an operad map $P\to\bar{\EE}_X$. If we dispose only of a model structure for reduced operads, we tacitly assume that our objects are based, and that our endomorphism-operads are the reduced ones. The following theorem is a model-theoretic formulation of Boardman and Vogt's \emph{homotopy invariance property} \cite[thm. 4.58]{BV}.

\begin{thm}\label{invariance}Let $f:X\to Y$ be a (based) map between (based) objects of a monoidal model category in which the (reduced) operads carry a transferred model structure; let $P$ be a cofibrant operad.

(a) If $Y$ is fibrant and $f^{\otimes n}$ is a trivial cofibration for each $n\geq 1$, then any $P$-algebra structure on $X$ extends (along $f$) to a $P$-algebra structure on $Y$.

(b) If $X$ is cofibrant and $f$ is a trivial fibration, then any $P$-algebra structure on $Y$ lifts (along $f$) to a $P$-algebra structure on $X$.

(c) If $X$ and $Y$ are cofibrant-fibrant and $f$ is a weak equivalence, then any $P$-algebra structure on $X$ (resp. $Y$) induces a $P$-algebra structure on $Y$ (resp. $X$) in such a way that $f$ preserves the $P$-algebra structures up to homotopy.\end{thm} The latter statement means precisely that $f$ admits a factorization into a trivial cofibration $f_1:X\to Z$ followed by a trivial fibration $f_2:Z\to Y$ such that $f_1$ is a $P$-algebra map with respect to a structure map $\phi_1:P\to\EE_Z$, $f_2$ is a $P$-algebra map with respect to a structure map $\phi_2:P\to\EE_Z$, and the structure maps $\phi_1$ and $\phi_2$ are homotopic in the model category of operads (the homotopy relation is well defined since $P$ is a cofibrant operad and $\EE_Z$ is a fibrant operad).

\begin{proof}
We define a collection $\EE_{X,Y}$ by $\EE_{X,Y}(n)=Y^{(X^{\otimes n})}$. The endomorphism-operad $\EE_f$ of $f$ is defined by the following pullback-diagram of collections:\begin{diagram}[small]\EE_f&\rTo&\EE_X\\\dTo&&\dTo_{f_*}\\\EE_Y&\rTo^{f^*}&\EE_{X,Y}\end{diagram}(Set theoretically, $\EE_f(n)=\{(\phi,\psi)\in\EE_X(n)\times\EE_Y(n)\,|\,f\phi=\psi f^{\otimes n}\}$.) From the operads $\EE_X$ and $\EE_Y$ the collection $\EE_f$ inherits the structure of an operad. Moreover, $f$ is compatible with the $P$-algebra structure maps $P\to\EE_X$ and $P\to\EE_Y$ if and only if these are induced by an operad map $P\to\EE_f$. 

Since trivial fibrations are closed under pullback, the exponential transpose of the pushout-product axiom shows that under the hypothesis of (a), the horizontal maps of the above diagram are trivial fibrations. Therefore, since $P$ is a cofibrant operad, the $P$-algebra structure map $P\to\EE_X$ has a lift $P\to\EE_f\to\EE_Y$ which yields the required $P$-algebra structure on $Y$. Dually, the hypothesis of (b) implies that $\EE_f\to\EE_Y$ is a trivial fibration, whence the required lift of the $P$-algebra structure map $P\to\EE_Y$ to $P\to\EE_f\to\EE_X$.

Assume now that $f$ is a weak equivalence between cofibrant-fibrant objects and that $X$ comes equipped with a $P$-algebra structure. The weak equivalence $f$ factors into a trivial cofibration $f_1:X\to Z$ followed by a trivial fibration $f_2:Z\to Y$. Since $X$ and $Y$ are cofibrant, we may assume that $f_2$ admits a trivial cofibration as section; in particular, the tensor powers of $f_2$ are weak equivalences between cofibrant objects, cf. (\ref{stable}), and we get a pullback diagram of fibrant collections\begin{diagram}[small]\EE_{f_2}&\rTo^\phi&\EE_Z\\\dTo&&\dTo_{(f_2)_*}\\\EE_Y&\rTo^{(f_2)^*}&\EE_{Z,Y}\end{diagram}in which the vertical maps are trivial fibrations and the horizontal maps are weak equivalences. It follows that for the cofibrant operad $P$, the upper horizontal map $\phi$ induces a bijection between homotopy classes$$[P,\phi]:[P,\EE_{f_2}]\cong[P,\EE_Z]$$Since $f_1$ is a trivial cofibration, the $P$-algebra structure map $P\to\EE_X$ extends to a $P$-algebra structure map $\phi_1:P\to\EE_Z$. The latter has a (up to homotopy unique) lift $\psi:P\to\EE_{f_2}$ such that $\phi_1$ and $\phi_2=\phi\psi$ are homotopic. The composite map $P\to\EE_{f_2}\to\EE_Y$ yields the required $P$-algebra structure on $Y$.

A dual argument shows that a $P$-algebra structure on $Y$ induces a $P$-algebra structure on $X$ in such a way that $f$ preserves the $P$-algebra structures up to homotopy in the above mentioned sense.\end{proof}

\begin{rmk}\label{stable}The slight asymmetry between (\ref{invariance}a) and (\ref{invariance}b) is due to the fact that the tensor powers of a trivial cofibration are in general not trivial cofibrations, cf. (\ref{monoid}). The latter becomes true if the domain of the considered trivial cofibration is cofibrant, or more generally, if the monoidal model category has a generating set of trivial cofibrations with cofibrant domains, cf. (\ref{transfer}). Therefore, property (\ref{invariance}a) shows that if the generating trivial cofibrations of the underlying model category have cofibrant domains, then \emph{the category of algebras over a cofibrant operad admits a fibrant replacement functor}.\end{rmk}

\section{Comparison theorems}

Throughout this section, $\EE$ is a monoidal model category satisfying either the hypotheses of (\ref{operad}) or the hypotheses of (\ref{operad2}); in particular, $\EE$ is assumed to have a cofibrant unit. Thus, (reduced) operads in $\EE$ carry a model structure. An operad $P$ is called \emph{admissible} if the category of $P$-algebras carries a model structure which is transferred (\ref{transfer}) from $\EE$ along the free-forgetful adjunction $\FF_P:\EE\lrto\Alg_P:\UU_P$. Under mild assumptions on $\EE$, cf. (\ref{adm}), cofibrant operads are admissible. For an arbitrary operad $P$, we define ``the'' category of \emph{homotopy $P$-algebras} as the category of $\hat{P}$-algebras for some cofibrant replacement $\hat{P}$ of $P$. We will show that this category is well defined up to Quillen equivalence (Corollary \ref{invariant2} below). 

Recall that an operad $P$ is \emph{$\Sg$-cofibrant} if the collection underlying $P$ is cofibrant. The main purpose of this section is to show that for an \emph{admissible $\Sg$-cofibrant} operad $P$, the category of $P$-algebras and the category of homotopy $P$-algebras have equivalent homotopy categories. The class of admissible $\Sg$-cofibrant operads includes most of the commonly used $A_\infty$- and $E_\infty$-operads. The difference between cofibrant and admissible $\Sg$-cofibrant operads is reminiscent of the difference between projective and flat objects in homological algebra. 

An operad $P$ is called \emph{$\Sg$-split} if $P$ is retract of $P\otimes\Ass$ where $\Ass$ is the operad acting on associative monoids. A definition resembling this occurs in \cite{H}.\vspace{1ex}

We first give some criteria for an operad to be admissible, using again the path-object argument (\ref{path}). Observe in particular that (\ref{admissible}b) holds under the hypotheses of (\ref{operad}), and (\ref{admissible}c) holds under the hypotheses of (\ref{operad2}).

\begin{prp}\label{admissible}Let $\EE$ be a cofibrantly generated monoidal model category with cofibrant unit and symmetric monoidal fibrant replacement functor.\vspace{1ex} 

(a) If there exists an operad map $j:P\to P\otimes Q$ and an interval in $\EE$ with a $Q$-coalgebra structure, then $P$ is admissible.

(b) If there exists an interval in $\EE$ with a coassociative comultiplication, then $\Sigma$-split operads are admissible.

(c) If there exists an interval in $\EE$ with a coassociative and cocommutative comultiplication, then all operads are admissible.\end{prp}

\begin{proof} (a) implies (b) resp. (c), putting $Q=\Ass$ resp. $Q=\Com$, where $\Com(n)=I$ is the operad for commutative monoids and $\Ass(n)=I[\Sg_n]$ is the operad for associative monoids. 

For (a), we can use the path-object argument, since the forgetful functor $\Alg_P\to\EE$ preserves filtered colimits. Let $A$ be a $P$-algebra. The symmetric monoidal fibrant replacement functor $A\eqv \tilde{A}$ induces a $\tilde{P}$-structure on $\tilde{A}$, where $i:P\eqv \tilde{P}$ is the fibrant replacement for $P$. This yields a fibrant replacement functor $A\eqv i^*\tilde{A}$ for $P$-algebras. Moreover, for any $P$-algebra $A$, mapping the $Q$-coalgebraic interval $J$ into $A$ yields a $P\otimes Q$-algebra $A^J$ by (\ref{Hopf}). Thus, $j^*(A^J)$ defines a functorial path-object for fibrant $P$-algebras using (\ref{weak}), the pushout-product axiom and the hypothesis that the unit is cofibrant. Therefore, the model structure of $\EE$ transfers to $\Alg_P$ and $P$ is admissible.\end{proof}

\begin{rmk}\label{adm}Spitzweck \cite[thm. 4.3]{Sp} proves the admissibility of cofibrant operads under the hypothesis that the \emph{monoid axiom} of Schwede-Shipley \cite{SS} holds. In many examples, the admissibility of cofibrant operads may be established using (\ref{admissible}): all topological resp. simplicial operads are admissible by (\ref{admissible}c). Cofibrant chain operads are admissible by a construction of \cite[thm. 2.1.1]{mi}: indeed, in the category of $R$-chain operads there exists a $\Sigma$-cofibrant resolution $\EE_\infty\overset{\sim}{\cto}\Com$ together with a canonical $\EE_\infty$-coalgebra structure on the standard $R$-chain interval. For each $R$-chain operad $P$, this induces a trivial fibration $P\otimes\EE_\infty\overset{\sim}{\cto}P$; for cofibrant operads $P$, the latter admits a section so that (\ref{admissible}a) applies.\end{rmk}



\begin{prp}\label{cof}Cofibrant operads are $\Sg$-cofibrant.\end{prp}

\begin{proof}Any cofibrant operad is retract of a cellular operad, i.e. a cellular extension of the initial operad, cf. (\ref{transfer}). Since $\Sg$-cofibrant operads are closed under retract, and since the initial operad is $\Sg$-cofibrant, it is enough to show that the class of $\Sg$-cofibrant operads is closed under cellular extensions: this is done in (\ref{free}).\end{proof}

\begin{thm}\label{invariant}In a left proper monoidal model category (with cofibrant unit), the base-change adjunction with respect to a weak equivalence of admissible $\Sg$-cofibrant operads is a Quillen equivalence.\end{thm}

\begin{proof}Let $\phi:P\to Q$ be a weak equivalence of admissible $\Sg$-cofibrant operads. The base-change adjunction $\phi_!:\Alg_P\lrto\Alg_Q:\phi^*$ is a Quillen pair since by the definition of the model structures, the restriction functor $\phi^*$ preserves weak equivalences and fibrations, so that its left adjoint $\phi_!$ preserves cofibrations. Since $\phi^*$ also reflects weak equivalences, the derived adjunction is an equivalence if (and only if) the unit induces a weak equivalence $A\to\phi^*\phi_!A$ for each \emph{cofibrant} $P$-algebra $A$. Since any cofibrant $P$-algebra is retract of a \emph{cellular $P$-algebra} (\ref{transfer}) and we assume that the model category is left proper, this follows from (\ref{cellular}).\end{proof}


We define a \emph{$\Sigma$-cofibrant resolution} of $P$ to be a $\Sigma$-cofibrant operad $P_\infty$ together with a trivial fibration of operads $P_\infty\overset{\sim}{\cto}P$. Recall that the category of homotopy $P$-algebras is the category of $\hat{P}$-algebras for some cofibrant replacement $\hat{P}$ of $P$. 

\begin{cor}\label{invariant2}Assume that cofibrant operads are admissible and that the underlying model category is left proper (and has a cofibrant unit). Then for any admissible $\Sg$-cofibrant resolution $P_\infty$ of $P$, the category of $P_\infty$-algebras is Quillen equivalent to the category of homotopy $P$-algebras.\end{cor}

\begin{proof}A trivial fibration $P_\infty\overset{\sim}{\cto}P$ induces a weak equivalence $\hat{P}\eqv P_\infty$ for any cofibrant replacement $\hat{P}$ of $P$. Since $P_\infty$ is $\Sg$-cofibrant, (\ref{cof}) and (\ref{invariant}) imply that the category of $\hat{P}$-algebras is Quillen equivalent to the category of $P_\infty$-algebras.\end{proof}

\begin{rmk}An admissible $\Sg$-cofibrant resolution of the operad $\Ass$ (resp. $\Com$) is a so called \emph{$A_\infty$-operad} (resp. \emph{$E_\infty$-operad}). The corresponding algebras are $A_\infty$-algebras (resp. $E_\infty$-algebras). Under the assumption of (\ref{invariant2}), the homotopy categories of $A_\infty$- resp. $E_\infty$-algebras are up to equivalence of categories independent of the chosen $A_\infty$- resp. $E_\infty$-operad.

Under the assumptions of (\ref{operad}) or (\ref{operad2}), the operad $\Ass$ is itself admissible $\Sg$-cofibrant. Indeed, the underlying collection is cofibrant since the unit of $\EE$ is cofibrant. Moreover, there is a diagonal $\Ass\to \Ass\otimes\Ass$ so that (\ref{admissible}a) implies that $\Ass$ is admissible. In other words, in any left proper monoidal model category satisfying our hypotheses, \emph{associative monoids carry a transferred model structure} and by (\ref{invariant}) \emph{$A_\infty$-algebras may be rectified to monoids without loss of homotopical information}. In the topological case, this has been established by Stasheff \cite{St}, Boardman and Vogt \cite[thm. 1.26]{BV} and May \cite[thm. 13.5]{M}. Schwede and Shipley \cite[thm. 3.1]{SS} prove the existence of a transferred model structure for associative monoids under the assumption that the \emph{monoid axiom} holds.\vspace{1ex}

The homotopy theory of $E_\infty$-algebras is more involved, since the operad $\Com$ is not $\Sg$-cofibrant. We discuss examples (\ref{examples}): 

\begin{stit}\emph{Simplicial $E_\infty$-algebras.}

The category of simplicial sets has a canonical $E_\infty$-operad given by the universal $\Sg_n$-bundles $\WW\Sg_n$ (not to be confused with Boardman and Vogt's $\WW$-construction). The operad structure is induced by the permutation operad, since the simplicial $\WW$-construction is product-preserving. The category of $\WW\Sg$-algebras has been extensively studied by Barratt and Eccles \cite{BE} for the construction of their infinite loop space machine. The comparison theorem (\ref{comparison}) below relates their simplicial approach to the more classical topological approach.\end{stit}

\begin{stit}\emph{Topological $E_\infty$-algebras.}

The geometric realization of $\WW\Sg$ is a topological $E_\infty$-operad. Boardman and Vogt's \emph{little cubes operad} $\CC_\infty$ is unlikely to be $\Sg$-cofibrant for the model structure we consider, because of the lack of a suitable equivariant $CW$-structure. Nonetheless, since the $\Sg$-actions are free, this operad has similar invariance properties as $\Sg$-cofibrant operads, cf. Vogt \cite{V}. The importance of the little cubes operad stems from its canonical action on \emph{infinite loop spaces}; the latter \emph{fully} embed in $\CC_\infty$-algebras and are characterized up to homotopy as the group-complete $\CC_\infty$-algebras.\end{stit}

\begin{stit}\emph{$E_\infty$-chain algebras.}

Since the normalized $R$-chain functor $\NN^R_*$ is symmetric monoidal, it sends operads to operads. The normalized $R$-chains $\EE_\infty=\NN^R_*(\WW\Sg)$ form thus a $\Sg$-free resolution of $\Com$. Since $\NN^R_*$ is also comonoidal, we get a diagonal $\EE_\infty\to\EE_\infty\otimes\EE_\infty$. Fresse and the first named author construct in \cite[thm. 2.1.1]{mi} an $\EE_\infty$-coalgebra structure on the chains (or dually, a $\EE_\infty$-algebra structure on the cochains) of any simplicial set. It follows from (\ref{admissible}a), that $\EE_\infty$ is admissible, i.e. an $E_\infty$-operad for the category of $R$-chain complexes. According to Mandell \cite[main thm.]{Man}, an $E_\infty$-structure on the cochains of a nilpotent simplicial set $X$ is a complete invariant of the \emph{$p$-adic homotopy type} of $X$, provided that $R$ is a field of characteristic $p$ with surjective Frobenius map.\end{stit}

\begin{stit}\emph{$E_\infty$-ring spectra.} 

No monoidal model category for \emph{stable homotopy} can simultaneously satisfy the first two hypotheses of (\ref{operad}), by a well known argument due to Lewis, cf. \cite[XIV]{MMSS}. However, all known models for stable homotopy are enriched either in simplicial sets or in topological spaces. Therefore, it makes sense to speak of simplicial (or topological) operad actions on spectra. In the enriched case, Quillen's axiom $SM7$ \cite[II.2]{Q} replaces the pushout-product axiom and guarantees that Theorem (\ref{invariance}) remains true for a cofibrant simplicial (or topological) operad $P$ and a general map of spectra $f:X\to Y$.

Property (\ref{invariance}a) implies the existence of a fibrant replacement functor for the category of spectra with $P$-algebra structure, provided that the generating trivial cofibrations of the monoidal model category have cofibrant domains, cf. (\ref{stable}). The suspension spectrum functor endows the category of spectra with an interval with coassociative and cocommutative comultiplication. The argument of (\ref{admissible}c) then yields a model structure on the category of spectra with $P$-algebra structure. In particular, there is a model structure on $A_\infty$- resp. $E_\infty$-ring spectra in any of the considered models for stable homotopy, provided that the chosen simplicial (or topological) $A_\infty$- resp. $E_\infty$-operad is cofibrant. Also, by the same argument as above, any $A_\infty$-ring spectrum may be rectified to an associative ring spectrum.\end{stit}\end{rmk}\vspace{1ex}

We conclude this section with the following comparison theorem for algebras in Quillen equivalent model categories. For the precise statement we need the following definition: A \emph{monoidal Quillen equivalence} between monoidal model categories (\ref{monoid}) is a Quillen pair, which is simultaneously a Quillen equivalence and a monoidal adjunction. An adjunction between symmetric monoidal categories is said to be monoidal if the left and right adjoint functors are symmetric monoidal, cf. (\ref{monoidal}), and moreover the unit and counit of the adjunction are monoidal transformations. Observe that if the left adjoint of an adjoint pair is a strong symmetric monoidal functor (\ref{monoidal}), then the right adjoint carries a natural symmetric monoidal structure for which the adjunction is monoidal.

\begin{thm}\label{comparison}Let $(L,R)$ be a monoidal Quillen equivalence between monoidal model categories in which the (reduced) operads carry a model structure. Let $P$ be a (reduced)  operad in the domain of $L$. Assume either that $L$ preserves weak equivalences or that $P(n)$ is cofibrant for all $n$. Then, the homotopy categories of homotopy $P$-algebras and of homotopy $LP$-algebras are equivalent.\end{thm}

\begin{proof}Since $L$ is symmetric monoidal, $L$ maps operads to operads. Let $\hat{P}$ be a cofibrant replacement of $P$. It follows by adjunction that $L\hat{P}$ is a cofibrant operad, which by either of the two hypotheses is a cofibrant replacement for $LP$. In particular, both $\hat{P}$ and $L\hat{P}$ are admissible. The functor $L$ maps $\hat{P}$-algebras to $L\hat{P}$-algebras. The functor $R$ maps  $L\hat{P}$-algebras to $RL\hat{P}$-algebras, which we consider as $\hat{P}$-algebras via the unit $\hat{P}\to RL\hat{P}$. Therefore, the adjunction $(L,R)$ lifts to an adjunction between $\hat{P}$-algebras and $L\hat{P}$-algebras. 

The given Quillen equivalence $(L,R)$ has the characteristic property that for cofibrant objects $A$ and fibrant objects $B$, a map $LA\to B$ is a weak equivalence if and only if the adjoint map $A\to RB$ is, cf. \cite[prp. 1.3.13]{Hov}. Assume now that $A$ is a cofibrant $\hat{P}$-algebra and $B$ a fibrant $L\hat{P}$-algebra. Since, according to (\ref{cof}), $\hat{P}$ has an underlying cofibrant collection, it follows from (\ref{cof2}) that $A$ has an underlying cofibrant object; moreover, $B$ has an underlying fibrant object. Therefore, a $L\hat{P}$-algebra map $LA\to B$ is a weak equivalence if and only if the adjoint $\hat{P}$-algebra map $A\to RB$ is, which establishes the equivalence of the homotopy categories of homotopy $P$-algebras and homotopy $LP$-algebras.\end{proof}

The preceding theorem shows in particular that the notion of an $E_\infty$-algebra is invariant under a monoidal Quillen equivalent change of the base category. In particular, the homotopy category of simplicial $E_\infty$-algebras is equivalent to the homotopy category of topological $E_\infty$-algebras.\newpage

\section{Appendix}

This Appendix establishes several key properties of $\Sg$-cofibrant operads which are used in Section $4$. These properties may to a large extent be derived from Spitzweck's work \cite[I.3.5 and I.4.5]{Sp}. Since his treatment uses the language of \emph{semi-model structures}, a topic we have not treated in this article, we give self-contained proofs of those model-theoretic properties we need. The hard work is actually concentrated in the proofs of (\ref{free0}) and (\ref{relative}), which we defer to the end of this Appendix.

A map of operads $\phi$ is called a \emph{$\Sg$-cofibration}, if $\UU(\phi)$ is a cofibration, where $\UU:\Oper(\EE)\to\Coll(\EE)$ is the forgetful functor from operads to collections. Recall that an operad $P$ is called $\Sg$-cofibrant, if $\UU(P)$ is cofibrant. A map of operads is called a \emph{free cofibration} if it is the image of a cofibration of collections under the free functor $\FF:\Coll(\EE)\to\Oper(\EE)$. Any pushout of a free cofibration is called a \emph{cellular extension}, cf. (\ref{transfer}). The reader may observe that the concept of a cellular extension makes sense even if there is no model structure on operads. 

Given an operad $P$ and a cofibration of collections $u:\UU(P)\to \KK$, we introduce a special notation for the pushout of $\FF(u):\FF\UU(P)\to\FF(\KK)$ along the counit $\FF\UU(P)\to P$, namely $P\to P[u]$, emphasizing that this pushout represents the cellular extension of $P$ generated by $u$. 

Similarly, if $A$ is a $P$-algebra and $u:\UU_P(A)\to Z$ is a cofibration, we denote by $A[u]$ the cellular $P$-algebra extension of $A$ generated by $u$, i.e. the pushout of $\FF_P(u):\FF_P\UU_P(A)\to\FF_P(Z)$ along the counit $\FF_P\UU_P(A)\to A$ of the adjunction $\FF_P:\EE\lrto\Alg_P:\UU_P$.

\begin{prp}\label{free0}For any $\Sg$-cofibrant operad $P$ and any cofibration of collections $u:\UU(P)\to\KK$, the cellular extension $P\to P[u]$ is a $\Sg$-cofibration.\end{prp}We postpone the proof to (\ref{proof1}).

\begin{cor}\label{free}A cellular extension of operads with $\Sg$-cofibrant domain is a $\Sg$-cofibration. The class of $\Sg$-cofibrant operads is thus closed under cellular extension.\end{cor}

\begin{proof}Given a cofibration of collections $\LL_1\to\LL_2$ and a map $\FF(\LL_1)\to P$, the cellular extension $P\to P\cup_{\FF(\LL_1)}\FF(\LL_2)$ may be identified with the cellular extension $P\to P[u]$ with respect to $u:\UU(P)\to\UU(P)\cup_{\LL_1}\LL_2$. The latter map is a pushout in collections along the adjoint $\LL_1\to\UU(P)$ of the given map $\FF(\LL_1)\to P$. Since the category of collections carries a model structure, $u$ is a cofibration so that (\ref{free0}) proves the assertion.\end{proof}

For an operad $P$ and a cofibration $u(0):P(0)\to Z$, define a collection $P_Z$ by $P_Z(0)=Z$ and $P_Z(n)=P(n)$ for $n>0$ and extend $u(0)$ to a cofibration of collections $u:\UU(P)\to P_Z$ setting $u(n)=id_{P(n)}$ for $n>0$.

\begin{lma}\label{under0}For an operad $P$ endowed with a cofibration $u(0):P(0)\to Z$, the category of $P[u]$-algebras is equivalent to the category of $P$-algebras under the cellular $P$-algebra extension $P(0)[u(0)]$.\end{lma}
\begin{proof}Recall that $P(0)$ comes equipped with a natural $P$-algebra structure making it the initial $P$-algebra. The asserted equivalence of categories breaks into two equivalences$$\Alg_{P[u]}\sim u/\Alg_{P}\sim P(0)[u(0)]/\Alg_P$$The category $u/\Alg_P$ has as objects the pairs $(A,v)$ consisting of a $P$-algebra $A$ and a map $v:Z\to\UU_P(A)$ such that $v\circ u(0)$ underlies the unique $A$-algebra map $P(0)\to A$. Morphisms $(A,v)\to(B,w)$ in $u/\Alg_P$ are given by maps $f:A\to B$ such that $w=\UU_P(f)v$. The first equivalence above is induced by pulling back the $P[u]$-algebra structure along the canonical operad maps $P\to P[u]$ and $\FF P_Z\to P[u]$. The second equivalence follows merely from the definition of the $P$-algebra extension $P(0)[u(0)]$.\end{proof}

\begin{prp}\label{under}Let $P$ be a $\Sg$-cofibrant operad and $A$ be a cellular $P$-algebra. Then there exists an operad $P[A]$ and a $\Sg$-cofibration of operads $\phi_A:P\to P[A]$ such that the category of $P$-algebras under $A$ is equivalent to the category of $P[A]$-algebras, and such that the following diagram commutes:\begin{diagram}[small]A/\Alg_P&\rTo^\sim&\Alg_{P[A]}\\\dTo^{U_A}&\ldTo_{\phi_A^*}&\\\Alg_{P}&&&\end{diagram}\end{prp}

\begin{proof} The proof is by induction on $A$, and naturally falls apart into three steps:\vspace{1ex}

\noindent(i) If $A$ is the initial $P$-algebra $P(0)$, we define $P[A]=P$ and $\phi_A=id_P$.\vspace{1ex}

\noindent(ii) Assume that $B$ is a cellular $P$-algebra extension of $A$, constructed as a pushout of $\FF_P(X)\to \FF_P(Y)$ along $\FF_P(X)\to A$ in $\Alg_P$, where $X\ito Y$ is a generating cofibration of the underlying model category $\EE$. Assume by induction that we have already constructed a $\Sg$-cofibration $\phi_A:P\to P[A]$ with the required properties. 

Notice first that the $P$-algebra map $A\to B$ makes $B$ into a $P[A]$-algebra. In particular, the category of $P$-algebras under $B$ is equivalent to the category of $P[A]$-algebras under $B$. Moreover, $B$ can be constructed as the cellular $P[A]$-algebra extension $A\to A[u(0)]$ where $u(0):\UU_{P[A]}(A)\to\UU_{P[A]}(A)\cup_X Y$ is induced by the adjoint of $\FF_{P[A]}(X)\to A$.

Since $P[A](0)=A$, it follows from (\ref{under0}) that the category of $P[A]$-algebras under $B$ is equivalent to the category of $P[A][u]$-algebras. We therefore let $P[B]=P[A][u]$. By (\ref{free0}), the canonical map $\phi_{BA}:P[A]\to P[B]$ is a $\Sg$-cofibration. We let $\phi_B$ be the composite map $\phi_{BA}\phi_A$.\vspace{1ex} 

\noindent(iii) Assume that $A$ is the colimit of a sequence of cellular $P$-algebra extensions$$A_0\ito A_1\ito A_2\ito\cdots$$indexed by some limit ordinal $\lambda$, and that we have constructed $\Sg$-cofibrations$$P[A_0]\to P[A_1]\to P[A_2]\to\cdots$$which identify $P[A_\xi]$-algebras with $P$-algebras under $A_\xi$ in a compatible way. Then we can simply take $P[A]$ to be the colimit of this sequence of operads, $P[A]=\mathrm{colim}_{\xi<\lambda}P[A_\xi]$. The canonical map $\phi_A:P\to P[A]$ is a $\Sg$-cofibration since the forgetful functor from operad to collections preserves filtered colimits.\end{proof}

\begin{cor}\label{cof2}Any cofibrant algebra over an admissible $\Sg$-cofibrant operad has a cofibrant underlying object.\end{cor}

\begin{proof}Since over admissible operads $P$, any cofibrant algebra is a retract of a cellular $P$-algebra $A$, it suffices to prove that the latter has a cofibrant underlying object. But, according to (\ref{under}), $A=P[A](0)$ and $P[A]$ is $\Sg$-cofibrant.\end{proof}

\begin{prp}\label{relative}Let $\phi:P\to Q$ be a weak equivalence of $\Sg$-cofibrant operads, let $u:\UU(P)\to\KK,\,v:\UU(Q)\to\LL$ be cofibrations and let $\psi:\KK\to\LL$ be a weak equivalence making the following square commutative:\begin{diagram}[small]\UU(P)&\rTo^{\UU(\phi)}&\UU(Q)\\\dTo^u&&\dTo_v\\\KK&\rTo^\psi&\LL\end{diagram}Then the induced map $P[u]\to Q[v]$ is a weak equivalence of operads.\end{prp}We postpone the proof to (\ref{proof2}).

\begin{prp}\label{cellular}In a left proper monoidal model category, the unit of the base-change adjunction with respect to a weak equivalence of $\Sg$-cofibrant operads $\phi:P\to Q$ induces a weak equivalence $A\to\phi^*\phi_!A$ at each cellular $P$-algebra $A$.\end{prp}

\begin{proof} We shall inductively construct weak equivalences of $\Sg$-cofibrant operads $\phi[A]:P[A]\to Q[\phi_!A]$ such that the map underlying the unit $A\to\phi^*\phi_!A$ may be identified with $\phi[A](0)$. As in (\ref{under}), the proof falls apart into three steps:\vspace{1ex}

\noindent(i) If $A$ is the initial $P$-algebra $P(0)$, then $\phi[A]=\phi$.\vspace{1ex}

\noindent(ii) Inductively, we assume given a weak equivalence $\phi[A]:P[A]\to Q[\phi_!A]$. We have to consider the case where $B$ is a cellular $P$-algebra extension $A\cup_{\FF_P(X)}\FF_P(Y)$ like in the proof of (\ref{under}ii). This gives rise to the following pushout diagrams, where $X\ito Y$ is a generating cofibration of the underlying model category $\EE$:\begin{diagram}[small]X&\rTo&P[A](0)&\rTo^{\phi[A](0)}&Q[\phi_!A](0)\\\dTo&&\dTo_{u(0)}&&\dTo_{v(0)}\\Y&\rTo& Z_A&\rTo^{\psi(0)}& Z_{\phi_!A}\end{diagram}The right pushout square induces by the definition preceding (\ref{under0}) a commutative square of collections\begin{diagram}[small]\UU(P[A])&\rTo^{\UU(\phi[A])}&\UU(Q[\phi_!A])\\\dTo_u&&\dTo_v\\P[A]_{Z_A}&\rTo^\psi&Q[\phi_!A]_{Z_{\phi_!A}}\end{diagram}The induction hypothesis implies that $\phi[A]$ is a weak equivalence; the left properness of $\EE$ implies that $\psi(0)$ (and hence $\psi$) is a weak equivalence. It follows from (\ref{relative}) that the induced map of operads $P[A][u]\to Q[\phi_!A][v]$ is a weak equivalence. By (\ref{under0}), the latter map may be identified with $\phi[B]:P[B]\to Q[\phi_!B]$.\vspace{1ex}

\noindent(iii) The case of sequential colimits simply comes down to the fact that if\begin{diagram}[small]P[A_0]&\rTo&P[A_1]&\rTo& P[A_2]&\rTo\cdots\\\dTo&&\dTo&&\dTo\\Q[\phi_!A_0]&\rTo&Q[\phi_!A_1]&\rTo&Q[\phi_!A_2]&\rTo\cdots\end{diagram}is a ladder with horizontal $\Sg$-cofibrations and vertical weak equivalences, then the map induced on the colimit operads $\mathrm{colim}P[A_\xi]\to\mathrm{colim}Q[\phi_!A_\xi]$ is again a weak equivalence, because the forgetful functor from operads to collections preserves filtered colimits and all operads of the ladder are $\Sg$-cofibrant, cf. (\ref{model}).\end{proof}

\begin{tit}\emph{The free operad generated by a collection.}\label{graft}\end{tit}

We discuss some preliminaries on the free functor $\FF:\Coll(\EE)\to\Oper(\EE)$. The free functor is based on an operation known as \emph{grafting of trees}, as explained e.g. in \cite{GK}, \cite[I.1]{GiK} or \cite[def. 1.37]{MSS}. To describe this grafting operation, we introduce the \emph{groupoid} $\TT$ of planar trees and non-planar isomorphisms: 

The \emph{objects} of $\TT$ are finite rooted planar trees. We closely follow the convention of Getzler-Kapranov \cite{GK}. Each edge in the tree has a natural orientation, so that we can speak of a vertex being at the beginning or at the end of an edge. Any tree will have three kinds of edges, namely \emph{internal edges} with a vertex at the beginning as well as at the end of the edge, \emph{input edges} with a vertex only at the end, and one outgoing edge, called the \emph{output} of the tree, with the root vertex as its beginning and no vertex at its end. The input edges and the output edge are together referred to as the \emph{external edges} of the tree. In addition, we will also need the tree pictured $|$, with no vertex and just one edge which is at the same time input and output; this tree serves as unit for the grafting operation on trees. The number of edges ending in a given vertex $v$ is called the \emph{valence} of $v$ and denoted $|v|$. A vertex with valence $0$ is called a stump. The set of input edges of the tree $T$ is denoted $in(T)$; this set has a natural linear ordering inherited from the planar structure of the tree. The cardinality of $in(T)$ is denoted $|T|$. 

Here is a picture of a tree with $4$ vertices, $3$ internal edges and $5$ input edges. The root has valence $2$; there are two vertices of valence $3$ and a stump.

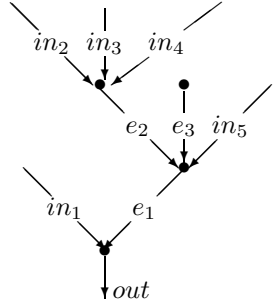
\begin{diagram}[height=0.55cm,abut]&&&&&\\\rdTo(1.2,2)~{{in}_2}\dTo~{{in}_3}\hspace{0.5em}\ldTo(1.3,2)~{{in}_4}&&\\\bullet\,\,&\bullet&&\\&\rdTo(1.2,2)~{e_2}\dTo~{e_3}\ldTo(1.2,2)~{{in}_5}&\\&\bullet&&&\\\rdTo(1,2)~{{in}_1}\ldTo(1,2)~{e_1}&&&\\\bullet&&&\\\dTo_{out}&&&&\end{diagram}

The \emph{morphisms} of $\TT$ are isomorphisms of trees, where we forget the planar structure. In particular, any isomorphism $\phi:T\to T'$ maps vertices to vertices, the root to the root, internal edges to internal edges and inputs to inputs, thus $|T|=|T'|$ and $|v|=|\phi(v)|$.

Any tree $T$ with a root of valence $n$ decomposes canonically into $n$ trees $T_1,\dots,T_n$ whose outputs are grafted upon the inputs of the tree $t_n$ with one vertex and $n$ inputs. We denote this grafting operation by $T=t_n(T_1,\dots,T_n)$. Observe that the number of vertices of each $T_i$ is strictly less than the number of vertices of $T$, which allows for inductive definitions. Any isomorphism $\phi:T\to T'$ decomposes as $\phi=\sg(\phi_1,\dots,\phi_n)$ with isomorphisms $\sg:t_n\to t_n$ and $\phi_i:T_i\to T'_{\sg(i)},\,i=1,\dots,n$. We identify the automorphism group of $t_n$ with the symmetric group $\Sg_n$.

For any collection $\KK$  we define a contravariant functor $\underline{\KK}:\TT^\op\to\EE$ putting inductively $\underline{\KK}(|)=I$ (the unit of the underlying monoidal model category $\EE$) and $$\underline{\KK}(T)=\underline{\KK}(t_n(T_1,\dots,T_n))=\KK(n)\otimes\underline{\KK}(T_1)\otimes\cdots\otimes\underline{\KK}(T_n).$$
On morphisms $\phi:T\to T'$, we get again by induction$$\phi^*=\sg(\phi_1,\dots,\phi_n)^*=\sg^*\otimes\phi^*_{\sg^{-1}(1)}\otimes\cdots\otimes\phi^*_{\sg^{-1}(n)}.$$There is also a covariant set-valued functor $\lambda:\TT\to\Sets$ associating to each tree $T$, the set $\lambda(T)$ of numberings of $\inp(T)$. A numbering $\tau\in\lambda(T)$ consists of a bijection $\tau:\{1,\dots,|T|\}\to\inp(T)$. Any isomorphism $\phi:T\to T'$ induces (by composition with $\inp(\phi):\inp(T)\to\inp(T')$) a bijection $\lambda(\phi):\lambda(T)\to\lambda(T')$. Since the category of sets naturally maps to $\EE$ via $S\mapsto\coprod_{s\in S}I$, we can consider $\lambda$ as a covariant functor with values in $\EE$. The classical formula for the free operad $\FF\KK$ generated by the collection $\KK$ amounts to the following tensor product over the groupoid $\TT$: $\FF\KK=\underline{\KK}\otimes_\TT\lambda.$
Since $\TT$ falls apart as a disjoint sum of groupoids $\TT(n)=\{T\in\TT\,|\,|T|=n\}$, $\FF\KK$ is the sum of the objects $\FF\KK(n)=\underline{\KK}\otimes_{\TT(n)}\lambda,\,n\geq 0$; in particular, it has a natural structure of collection where the symmetric groups $\Sg_{|T|}$ act from the right on the numberings $\tau\in\lambda(T)$. This restricted tensor product $\underline{\KK}\otimes_{\TT(n)}\lambda$ again decomposes as a sum, indexed by isomorphism classes of trees in $\TT(n)$:$$\FF\KK(n)=\coprod_{[T]\in\TT(n)/\sim}\underline{\KK}(T)\otimes_{\Aut(T)}I[\Sg_n]$$By categorical generalities (using the isomorphism $X\otimes I[\Sg_n]\cong\coprod_{\sg\in\Sg_n}X$), the tensor product may also be identified with the colimit$$\FF\KK=\mathrm{colim}_{\TT[\lambda]}\underline{\KK}\pi$$where $\pi:\TT[\lambda]\to\TT$ is the Grothendieck-construction applied to $\lambda:\TT\to\Sets$. Explicitly, $\TT[\lambda]$ is a groupoid whose objects are pairs $(T,\tau)$ consisting of a tree and a numbering of its inputs, and whose morphisms $\phi:(T,\tau)\to(T',\tau')$ are isomorphisms $\phi:T\to T'$ such that $\lambda(\phi)\tau=\tau'$. This groupoid is again a disjoint sum of subgroupoids $\TT[\lambda](n)$. Grafting of trees according to the given numberings endows $\TT[\lambda]$ with the structure of an \emph{operad in groupoids}. It follows from the inductive definition of $\underline{\KK}$ that $\FF\KK$ inherits from $\TT[\lambda]$ a natural structure of operad in $\EE$, and it is well known that the operad thus defined is the free operad generated by the collection $\KK$, cf. \cite{GK} or \cite[1.9]{MSS}.

It is straightforward to deduce from the preceding discussion that for a cofibrant collection $\KK$, the free operad $\FF\KK$ is $\Sg$-cofibrant. Proposition (\ref{free0}) generalizes this fact, and we now prepare the proof thereof. We need one more concept which in this context goes back to Hinich \cite{H}, namely trees endowed with a distinguished subset of \emph{coloured} vertices. We represent such coloured trees by pairs $(T,c)$ consisting of a tree $T$ and subset $c$ of coloured (internal) vertices. This gives rise to the groupoid $\widehat{\TT}$ of coloured trees and isomorphisms preserving the colourings. A coloured tree $(T,c)$ is \emph{admissible} if any internal edge of $T$ has at least one coloured extremity. For a coloured tree $(T,c)$, the grafting operation $T=t_n(T_1,\dots,T_n)$ yields canonical colourings for $t_n$ and $T_i,\,i=1,\dots, n$.

Given a map $u:\KK_1\to\KK_2$ of collections, and a coloured tree $(T,c)$, we define inductively an object $\underline{u}(T,c)$ of the underlying monoidal model category $\EE$ by $$\underline{u}(T,c)=\begin{cases}\KK_1(n)\otimes\underline{u}(T_1,c_1)\otimes\cdots\otimes\underline{u}(T_n,c_n)\text{ if the root is uncoloured},\\\KK_2(n)\otimes\underline{u}(T_1,c_1)\otimes\cdots\otimes\underline{u}(T_n,c_n)\text{ if the root is coloured},\end{cases}$$where we have the grafting operation $T=t_n(T_1,\dots,T_n)$ with colourings $c_i$ of $T_i$ corresponding to the colouring $c$ of $T$.

For any inclusion $c'\subset c$ of vertex-sets of $T$, there is a canonical map $\underline{u}(T,c')\to\underline{u}(T,c)$ induced by $u:\KK_1\to\KK_2$. If $u$ is a cofibration between cofibrant collections, then $\underline{u}(T,c')\to\underline{u}(T,c)$ is a cofibration in $\EE$ by the pushout-product axiom, cf. (\ref{monoidal}). Moreover, the canonical ``latching'' map$$\mathrm{colim}_{c'\subsetneq c}\underline{u}(T,c')\overset{def}{=}\underline{u}^-(T,c)\to\underline{u}(T,c)$$ is a cofibration, since the colimit may be constructed as an iterated pushout. We actually need the following stronger result:

\begin{lma}\label{adjoin}For any cofibration of cofibrant collections $u:\KK_1\to\KK_2$, the latching map $\underline{u}^-(T,c)\to\underline{u}(T,c)$ is an $\Aut(T,c)$-cofibration.\end{lma}

\begin{proof}Observe first that the latching map is indeed a map of $\Aut(T,c)$-objects, where $\Aut(T,c)$ is the automorphism group of $(T,c)$ in $\widehat{\TT}$. We assume the grafting operation $T=t_n(T_1,\dots,T_n)$ with colourings $c_i$ of $T_i$ corresponding to the colouring $c$ of $T$. 

In order to determine $\Aut(T,c)$, the set $\{(T_1,c_1),\dots,(T_n,c_n)\}$ is partitioned into subsets of pairwise isomorphic coloured trees, say $\{(T^1_1,c^1_1),\dots,(T^1_{n_1},c^1_{n_1})\}\cup\cdots\cup\{(T^k_1,c^k_1),\dots,(T^k_{n_k},c^k_{n_k})\}$ with $n_1+\cdots+n_k=n$. It follows that\begin{align*}\Aut(T,c)&\cong(\Aut(T^1,c^1)^{n_1}\rtimes\Sg_{n_1})\times\cdots\times(\Aut(T^k,c^k)^{n_k}\rtimes\Sg_{n_k})\\&\cong(\Aut(T^1,c^1)^{n_1}\times\cdots\times\Aut(T^k,c^k)^{n_k})\rtimes(\Sg_{n_1}\times\cdots\times\Sg_{n_k})\\&\overset{def}{=}G\rtimes\Sg\end{align*}If the root of $T$ is uncoloured, the map $\underline{u}^-(T,c)\to\underline{u}(T,c)$ may be identified with $\KK_1(n)\otimes (A\to B)$ for a certain map of $G\rtimes\Sg$-objects, denoted $A\to B$, whose underlying map is a cofibration. If the root of $T$ is coloured, we get $\KK_1(n)\otimes B\cup_{\KK_1(n)\otimes A}\KK_2(n)\otimes A\to\KK_2(n)\otimes B$ for the same $A\to B$. In both cases, this yields $G\rtimes\Sg$-cofibrations by (\ref{tool}).\end{proof}

\begin{lma}\label{tool}Let $G,\Sg$ be finite groups with $\Sg$ acting from the right on $G$.\vspace{1ex}

For any $\Sg$-cofibration $X\ito Y$ and any map of right $G\rtimes\Sg$-objects $A\to B$ whose underlying map is a cofibration, the induced map$$(X\otimes B)\cup_{(X\otimes A)}(Y\otimes A)\to Y\otimes B$$ is a $G\rtimes\Sg$-cofibration, where $G\rtimes\Sg$ acts on $Y\otimes B$ by $(y\otimes b)^{(g,\sg)}=y^\sg\otimes b^{(g,\sg)}$.\end{lma}

\begin{proof}Using the adjunction between $G\rtimes\Sg$-equivariant maps $Y\otimes B\to Z$ and $\Sg$-equivariant maps $Y\to\Hom_G(B,Z)$, the statement of the lemma is (by exponential transpose) equivalent to the property that for any $G\rtimes\Sg$-equivariant trivial fibration $Z\to W$, the canonical $\Sg$-equivariant map$$\Hom_G(B,Z)\to\Hom_G(B,W)\times_{\Hom_G(A,W)}\Hom_G(A,Z)$$ is a trivial fibration. Therefore, we may assume without loss of generality that $\Sg$ is the trivial group. The statement of the lemma is then equivalent to the property that the canonical $G$-equivariant map $$Z^Y\to W^Y\times_{W^X}Z^X$$is a trivial fibration. This in turn follows from the pushout-product axiom, since under the hypothesis of the lemma, the $G$-equivariant map $$(X\otimes B)\cup_{(X\otimes A)}(Y\otimes A)\to Y\otimes B$$ has a cofibration as underlying map.\end{proof}

\tit\emph{Proof of Proposition (\ref{free0})}.\label{proof1}\vspace{1ex}

\noindent We construct $P\to P[u]$ as a sequential colimit of cofibrations of collections$$F_0\ito F_1\ito F_2\ito\cdots.$$Of course,$$F_0(n)=P(n),\,n\geq 0.$$The $\Sg_n$-object $F_k(n)$ is inductively defined by the pushout-diagram below, where $(T,c)$ ranges over the set $\AA_k(n)$ of isomorphism classes of \emph{admissible coloured trees with $n$ inputs and $k$ coloured vertices}. The vertical map on the left comes from the operad structure of $P$ and the inductive definition of $F_{k-1}(n)$:\begin{diagram}\coprod_{[T,c]\in\AA_k(n)}\underline{u}^-(T,c)\otimes_{\Aut(T,c)}I[\Sg_n]&\rTo&\coprod_{[T,c]\in\AA_k(n)}\underline{u}(T,c)\otimes_{\Aut(T,c)}I[\Sg_n]\\\dTo&&\dTo\\F_{k-1}(n)&\rTo&F_k(n)\end{diagram}By (\ref{adjoin}), the latching map $\underline{u}^-(T,c)\to\underline{u}(T,c)$ is a $\Aut(T,c)$-cofibration. The functor $-\otimes_{\Aut(T,c)}I[\Sg_n]$ is the left adjoint of a Quillen pair $\EE^{\Aut(T,c)}\lrto\EE^{\Sg_n}$ and preserves cofibrations. Therefore, the upper horizontal map is a $\Sg_n$-cofibration, and the induced map $F_{k-1}(n)\to F_k(n)$ is a $\Sg_n$-cofibration too. Since all objects of the sequence are $\Sg_n$-cofibrant, the sequential colimit $$P(n)\to P[u](n)\overset{def}{=}\mathrm{colim}_kF_k(n)$$ is a $\Sg_n$-cofibration. We thus get a cofibration of collections $P\to P[u]$. The operad structure on $P[u]$ is defined by grafting of coloured trees, using the operad structure of $P$ in order to get back tensor products over admissible trees. The required universal property of the operad map $P\to P[u]$ follows from its inductive construction: at each step we adjoin one more free operation labelled by an element of $\KK$.\qed

\tit\emph{Proof of Proposition (\ref{relative}).}\label{proof2}

According to the preceding proof, the induced operad map $P[u]\to Q[v]$ may be obtained as the sequential colimit of a ladder of maps of cofibrant collections \begin{diagram}[small]F_0&\rTo&F_1&\rTo& F_2&\rTo&\cdots\\\dTo&&\dTo&&\dTo\\G_0&\rTo&G_1&\rTo&G_2&\rTo&\cdots\end{diagram}with horizontal cofibrations. It suffices thus to show that the vertical maps of this ladder are weak equivalences, cf. (\ref{model}). By hypothesis, the left most vertical map is the given weak equivalence $P\to Q$. For each $k>0$ and each $n\geq 0$, the vertical component $F_k(n)\to G_k(n)$ is obtained as the pushout of the following diagram:\begin{diagram}F_{k-1}(n)&\lTo&\coprod_{[T,c]}\underline{u}^-(T,c)\otimes_{\Aut(T,c)}I[\Sg_n]&\rTo&\coprod_{[T,c]}\underline{u}(T,c)\otimes_{\Aut(T,c)}I[\Sg_n]&\\\dTo&&\dTo&&\dTo\\G_{k-1}(n)&\lTo&\coprod_{[T,c]}\underline{v}^-(T,c)\otimes_{\Aut(T,c)}I[\Sg_n]&\rTo&\coprod_{[T,c]}\underline{v}(T,c)\otimes_{\Aut(T,c)}I[\Sg_n]&\end{diagram}The vertical maps of the latter diagram are weak equivalences by the induction hypothesis and the pushout-product axiom; the two right horizontal maps are $\Sg_n$-cofibrations by (\ref{adjoin}). Since all objects of the diagram are $\Sg_n$-cofibrant, \emph{Reedy's patching lemma} \cite[prp. 12.11]{DHK} implies that the induced vertical component $F_k(n)\to G_k(n)$ is a weak equivalence too, and we are done.\qed

\vspace{5ex}

\noindent{\sc Universit\'e de Nice, Laboratoire J.-A. Dieudonn\'e, Parc Valrose, 06108 Nice Cedex, France.}\hspace{2em}\emph{E-mail:} cberger$@$math.unice.fr\vspace{2ex}

\noindent{\sc Mathematisch Instituut, Postbus 80.010, 3508 TA Utrecht, The Netherlands.}\hspace{2em}\emph{E-mail:} moerdijk$@$math.uu.nl

\end{document}